\documentclass[12pt]{article}

\usepackage[english]{babel}
\usepackage[usenames]{color}
\usepackage[cp1250]{inputenc}
\usepackage{amsfonts}
\usepackage{amsthm}
\usepackage{graphicx}
\usepackage{epsfig}

\theoremstyle{plain}
\newtheorem{df}{Definition}

\newtheorem{tw}[df]{Theorem}
\newtheorem{lem}[df]{Lemma}

\newtheorem{sps}[df]{Observation}
\newtheorem{ex}[df]{Example}

\newtheorem{hyp}[df]{Hypothesis}

\begin{document}
\newcommand{\bea}{\begin{eqnarray}}
\newcommand{\eea}{\end{eqnarray}}
\newcommand{\be}{\begin{equation}}
\newcommand{\ee}{\end{equation}}
\newcommand{\beas}{\begin{eqnarray*}}
\newcommand{\eeas}{\end{eqnarray*}}
\newcommand{\bs}{\backslash}
\newcommand{\bc}{\begin{center}}
\newcommand{\ec}{\end{center}}

\title{Complex base numeral systems.}

\author{Jarek Duda}

\date{\it \footnotesize Jagiellonian University, Reymonta 4, 30-059 Kraków, Poland, \\
\textit{email:} dudaj@interia.pl}

\maketitle

\begin{abstract}
In this paper will be introduced large, probably complete family of
complex base systems, which are 'proper' - for each point of the
space there is a representation which is unique for all but some zero
measure set. The condition defining this family is the periodicity -
we get periodic covering of the plane by fractals in hexagonal-type structure,
what can be used for example in image compression. There will be introduced full
methodology of analyzing and using this approach - both for the
integer part: periodic lattice and the fractional: attractor of
some IFS, for which the convex hull or properties like dimension of
the boundary can be found analytically. There will be also shown how
to generalize this approach to higher dimensions and found some
proper systems in dimension 3.
\end{abstract}

\section{Introduction}
Standard numeral systems allows us to represent points from
$\mathbb{R}$ in some base ($z$) as a sequence of digits:
$$x=\sum_i a_i z^i$$
In this paper will be shown how to generalize it to 2 dimensions, using the structure of
complex plane $\mathbb{C}=\mathbb{R}+j\mathbb{R}$. It was previously
done by Donald Knuth \cite{knu}, who has introduced imaginary base
systems, or in \cite{jam} where the $j-1$ base system was introduced
and where was shown that such representations
can be used to simplify operations involving complex numbers in today's microprocessors. \\
In this paper will be presented much larger family of complex base system,
which seems to be complete.\\

We expect from the numeral systems, that the representation function is surjective -
each point has at least one representation.\\
But we rather cannot expect injectivity - for example
$0.111111..._2=1.000000..._2$ - there is countable number of points
that have two different representations.\\
In the 2-dimensional case such set won't be countable - it will have
Hausdorf's dimension in $[1,2)$ range. I will assume only, that it's
Lebesgue's measure is 0.
\begin{df}
  \emph{We will call the numeral system} proper, \emph{if the representation function is surjective
 and}
pseudoijective:\\$\exists_{S\subset \mathbb{C}}:\ \mu_2(S)=0,
\forall_{x\in \mathbb{C}\backslash S}\ x$\emph{ has at most one representation}.
\end{df}
Choose some $\mathbb{N}\ni n\geq 2$.\\
We will use $\bar{n}:=\{0,..,n-1\} $ digits and some base $z\in
\mathbb{C}$.\\
We have to analyze two completely different sets:\\
\emph{Integral part}: $I_z=\{\sum_{i\geq 0} a_i z^i :a_i\in \bar{n},\
\exists_N\forall_{i>N}\ a_i=0\}$\\
\emph{Fractional part}: $F_z=\{\sum_{i<0} a_i z^i :a_i\in \bar{n}\}$\\
For the binary system $n=2,\ z=2,\ \bar{n}=\{0,1\},\ I_2=\mathbb{Z},\ F_2=[0,1]$.\\

\begin{figure}[h]
    \centering
        \includegraphics[width=15cm]{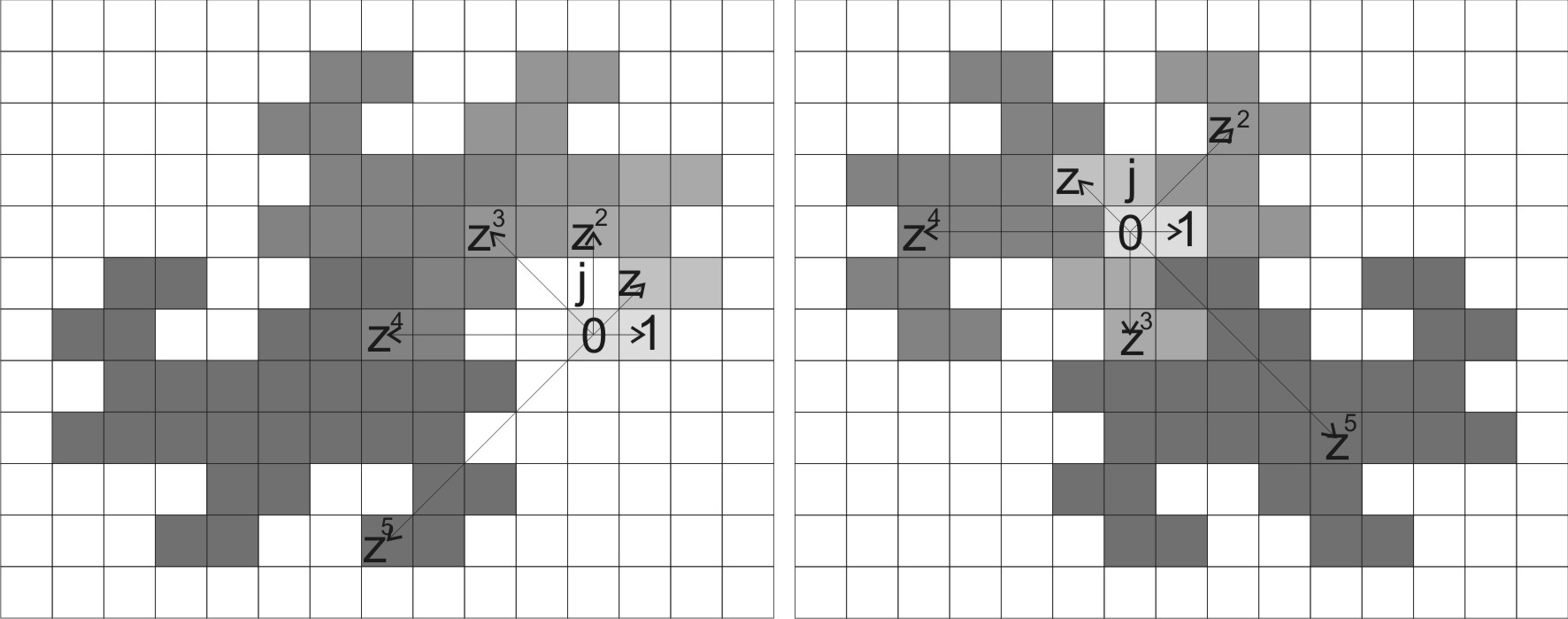}
        \caption{Integer part we can get using 6 youngest digits for $n=2$ and
        $z=j+1$ (left) or $z=j-1$ (right). It looks like in the first case we will
        make 'spirals' around 0, but in the second we should cover the whole space .}
\end{figure}

We will check, that they are fulfilling \emph{selfsimilarity
equations}: \be \label{ss1} z\cdot F_z=\bigcup_{i=0..n-1}
\left(F_z+i\right)\ee \be \label{ss2} \bigcup_{i=0..n-1}
\left(zI_z+i\right)=I_z.\ee For the binary system
$2[0,1]=[0,1]\cup([0,1]+1),\
2\mathbb{Z}\cup(2\mathbb{Z}+1)=\mathbb{Z}$.\\

\begin{figure}[h]
    \centering
        \includegraphics[width=15cm]{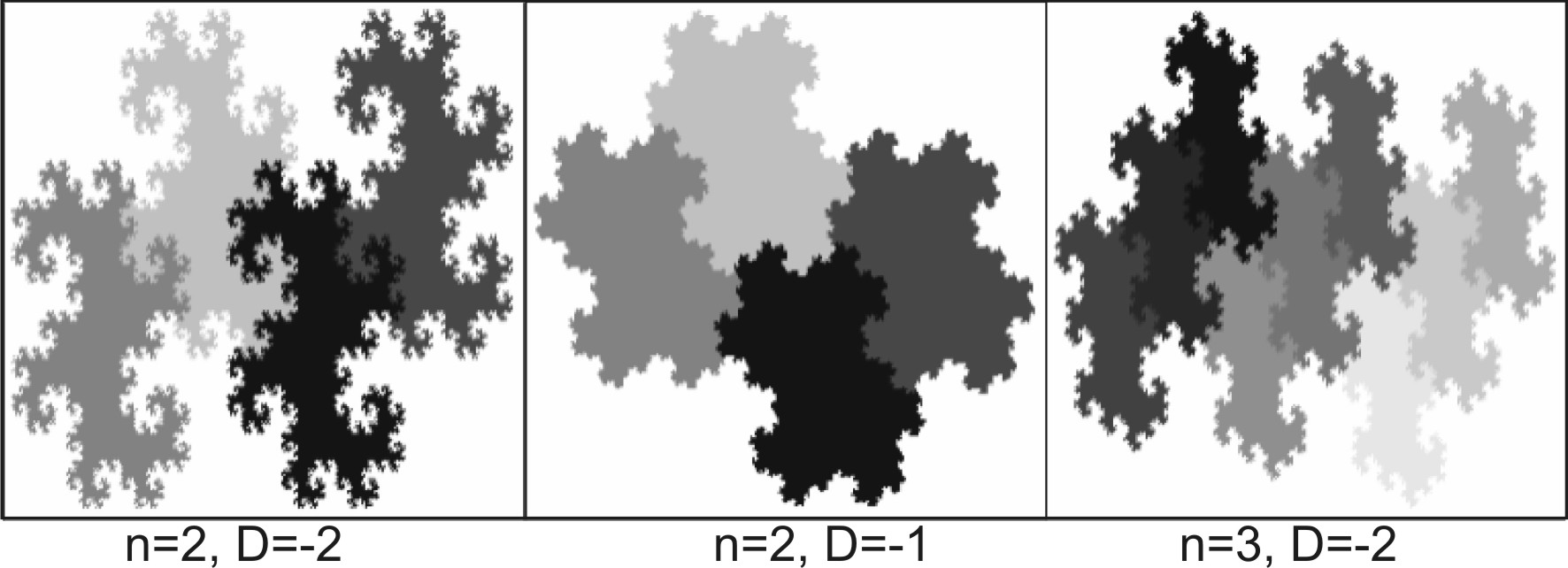}
        \caption{Some fractional parts. The color represents two oldest digits.} \label{ex}
\end{figure}

We will check that the necessity conditions implies:
\begin{lem} \label{nec}\emph{If} $(n,z)$ \emph{is proper}
$$\mu_2(F_z)>0,\quad \mu_2(F_z\cap (F_z+1))=0,$$
\be n=|z|^2. \label{powl} \ee
\end{lem}
\noindent It's why we will be dropping the upper index of $I_z^n,\ F_z^n$.\\

So $F_z$ is some compact, positive measure, central
symmetric(\cite{me}) set and when we take copies translated by all
points of some discrete set $I_z$ we 'tile' the whole plane, such
that tiles should intersects with its neighbors only on the boundaries.\\
Intuition says that it should be \emph{periodic tiling}:
\be \label{per} z^2\in \mathbb{Z}+z\mathbb{Z} \ee
We will see that it means, that $z=D/2+j\sqrt{n-(D/2)^2}$\\
for some $D\in \mathbb{Z}\cap \left(-2\sqrt n,2\sqrt n\right)$.

\begin{figure}[h]
    \centering
        \includegraphics[width=15cm]{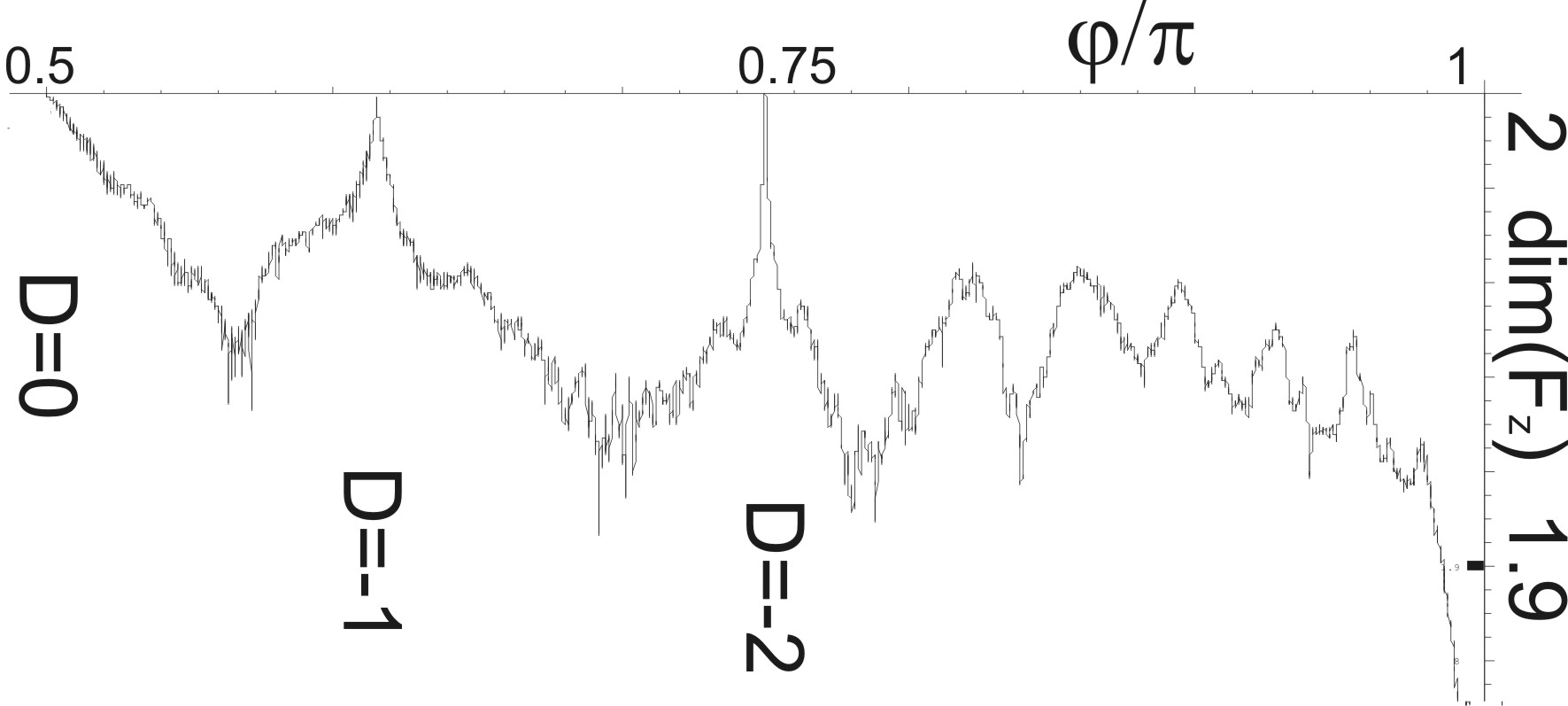}
        \caption{Numerically found Hausdorff's dimension of $F_z$.} \label{dimmm}
\end{figure}

\begin{hyp}
  All proper (2-dimensional) numeral systems are periodic (\ref{per}).
\end{hyp}
I couldn't prove it. On Fig. \ref{dimmm} are numerically found Hausdorff's
dimensions of $F_{\sqrt 2e^{j\varphi}}$  for $\varphi\in [\pi/2,\pi]$.
The tree peaks with dimension two are the periodic cases, while
$\alpha\rightarrow \pi,\ F_z$ is becoming a segment - its dimension drops to 1.\\

We will concentrate on the periodic case in this paper, and show
that the picture looks like it seems (like on Fig. \ref{ex} or Escher's pictures) -
we have slanting two dimensional lattice
($X=\mathbb{Z}+z\mathbb{Z}$) and each of its parallelogram (build of
$1$ and $z$) corresponds to exactly one copy of $F_z$ (so
$\mu_2(F_z)=|z|\sin(\varphi)$), which is connected, simply
connected, $F_z=\overline{\textrm{int}(F_z)}$ and it intersects with
its 6 neighbors (for $\varphi\neq \pm \pi/2$)
on its boundary, which dimension can be calculated analytically.\\
These boundaries ($\delta F_z$) are the sets on which points has two
or three representations - we can define dense set on which we are
loosing injectivity - like for the binary system:
$$S=\bigcup_{k\in \mathbb{Z}} z^k\left(I_z+\delta F_z\right)$$
In one dimensional case, this set is countable - has dimension equal
0, it's less by 1 than space dimension.
The skewness makes that this time, we have a bit worse (see table on page 16).

The surjectiveness condition can be written in a few ways: \be
\mathbb{C}=F_z+I_z=\sum_{i>0}z^iF_z=\overline{\sum_{i<0}z^i I_z}
\label{sur}\ee we can check it:
\begin{itemize}
  \item by checking that $I_z=\mathbb{Z}+z\mathbb{Z}$
  \item by checking that $0\in \textrm{int}(F_z)$
\end{itemize}

In the next section will be introduced formalism and proven some basic
properties: selfsimilarity relations, $|z|^2=n$ condition, that it's enough to
analyze argument of $z$ in one quarter and be characterizedS the periodic cases.\\
In the third section we will focus on the integer part. We will show
in what cases $I_z=\mathbb{Z}+z\mathbb{Z}$, what is sufficient for
system to be proper and that in the rest of cases situation isn't
much worse. The main tool will be the reduction - the
analogue of shift right with throwing out the youngest digit.\\
In the forth section we will focus on the fractional part. Firstly
there will be presented results from \cite{me} - a methodology of
finding the convex hull of $F_z$ and its properties analytically. In
the second - main subsection, will be introduced the methodology of
constructing the approximation of the boundary of $F_z$, which
allows to show that the tiling looks 'nice', quickly approximate
$F_z$ or to calculate analytically the dimension of the boundary.\\
In the last section, will be briefly presented a
proposition of generalization of this approaches into higher dimensions and
shown a family of proper system in dimension 3.\\

We can use the complex base systems is image compression too - in hexagonal-type
structure, we should get better correlations. \\
The other thing is that calculating the transform of the whole image would cost too much, so
the image is usually split into a periodic lattice of squares, which
are processed separately. In this approach lossy compression
creates recognizable lattice of lines. Splitting into lattice of
fractal-like shapes shouldn't have this problem.\\
To encode inside such shape - $\{\sum_{0\leq i \leq N} a_i
z^i:a_i\in \bar{n}\}$, we can use wavelet transform as it was just
usual binary system, for example basic Haar wavelet for $n=2$ would
be: $f_k=2a_k-1$ if $\forall_{i>k}\ a_i=0$, 0 else.\\
The boundary of larger integer sets is created from boundaries of
smaller ones, so the wavelet functions vanishing on the boundaries
should behave well in this case.\\
]

\section{Basic definitions and properties}
\begin{df}\emph{Fix some} $2\leq n\in\mathbb{N},\ z\in \mathbb{C}: |z|>1$\\
Digits $\overline{n}:=\{0,1,..,n-1\}$,\\
Pseudorepresentatnions
$P:=\{(a_i)_{i\in \mathbb{Z}}:\forall_i a_i\in\mathbb{Z},\ \exists_N\forall_{i>N}\ a_i=0\}$,\\
Shift left(right): $\textrm{\emph{shr}}(a_i):=(a_{i-1}),\quad \textrm{\emph{shl}}(a_i):=(a_{i+1})$,\\
Representations $P^n:=\{(a_i)\in P:\forall_i a_i\in \overline{n}\}$,\\
$P^+:=\{(a_i)\in P:\forall_{i<0}\ a_i=0\},\quad P^{n+}:=P^+\cap P^n$,\\
$P^-:=\{(a_i)\in P:\forall_{i\geq 0}\ a_i=0\},\quad P^{n-}:=P^-\cap P^n$,\\
Representing function $\textrm{P}\ni (a_i)\rightarrow (a_i)_z:=\sum_i a_i z^i$,\\
Integer part $(I_z\equiv )I^n_z:=P^{n+}_z$,\\
Fractional part $(F_z\equiv )F^n_z=P^{n-}_z$.
\end{df}

We will use standard set arithmetic:\\
$A+B=\{a+b:a\in A,\ b\in B\},\ f(A)=\{f(a):a\in A\}$.
\begin{sps}:
$$\left(\textrm{\emph{shl}}(a_i) \right)_z=z\cdot(a_i)_z,\quad
\left(\textrm{\emph{shr}}(a_i) \right)_z=z^{-1}\cdot (a_i)_z$$
$$\textrm{\emph{shl}}\left(P^{n-} \right)=P^{n-}+(..,0,\bar{n},0,..), \quad
\textrm{\emph{shl}}\left(P^{n+} \right)+(..,0,\bar{n},0,..)=P^{n+}$$
\emph{So we have the }selfsimilary equations:
$$ z\cdot F^n_z=F^n_z+\overline{n},\quad z\cdot I^n_z+\overline{n}=I^n_z.$$
\end{sps}
\noindent We've just checked (\ref{ss1}) and (\ref{ss2}) conditions.

The surjectiveness condition means that $P^n_z=\mathbb{C}$
\be P^n_z=\sum_{k>0} z^k F_z=F_z+I_z \ee
it's the sum of countable many translated copies of $F_z$ so we must have $\mu_2(F_z)>0$.\\
Pseudoinjectivity implies that $\forall_{i=1..n-1}\mu_2(F_z\cap(F_z+i))=0$ - points in this intersection would
have at least two representations - with $0$ and with $i$ at zeroth position.\\
Multiplication by $z$ is rotation and rescaling. Rotation and translation doesn't change the Lebegue's measure,
and the scaling multiply it by the scaling factor to power dimension - from the selfsimilarity
condition we have:
$$n\cdot \mu_2(F_z)=\mu_2\left(\bigcup_{i=0..n-1} (F_z+i)\right)=\mu_2(z\cdot F_z)=|z|^2\mu_2(F_z)$$
So we've checked Lemma \ref{nec}.\\
Observe that if we assume $|z|^2=n$, from the derivation above:
$\mu_2(F_z)>0$ automatically implies pseudoinjectivity - we have only to check the surjectiveness.\\
For the rest of the paper denote: \be z=\sqrt n e^{j\varphi}.\ee

When we conjugate $z$, $F_z$ will be conjugated.\\
When we take $-z$, we get ($-\bar n=\bar n -(n-1)$)\\
$F_{-z}=\sum_{i<0}z^i(-1)^i\bar n=\sum_{i<0}z^i \bar n-\sum_{\textrm{\tiny{odd}}\ i<0}z^i(n-1)=
F_z-\frac{(n-1)z^{-1}}{1-z^{-2}}=F_z-\frac{(n-1)z}{z^2-1}:$
\begin{sps}
\be \label{con} F_{\bar z}=\bar{F_z}\ee
\be \label {neg} F_{-z}=F_z-\frac{(n-1)z}{z^2-1}. \ee
\end{sps}

We will now find possible $z$, fulfilling \emph{periodicity condition}
($z^2\in \mathbb{Z}+z\mathbb{Z}$) and $|z|^2=n$.\\
Thanks of (\ref{con}), we can restrict to positive imaginary part and eventually take
conjugation of everything.
$$z=d+j\sqrt{n-d^2}\quad \textrm{for some }d\in \mathbb{R}$$
$$z^2=d^2-n+d^2+2jd\sqrt{n-d^2}$$
Comparing imaginary parts, we have $D:=2d\in \mathbb{Z}$, we finally have:
\begin{sps} \label{sps1}
  All periodic cases are:
\be z=D/2+j\sqrt{n-(D/2)^2} \ee
\be z^2=Dz-n \ee
for some $2\leq n\in \mathbb{N},\ D\in \mathbb{Z}\cap \left(-2\sqrt n,2\sqrt n\right).$
\end{sps}

In the next section we will prove:
\begin{tw}\label{main}
  Complex system fulfilling the periodicity condition is proper iff
  $ D\leq 1$.
\end{tw}

\section{The integral part}
In this part we will check whether $I_z=X$, where we will denote
$$X:=\mathbb{Z}+z\mathbb{Z}.$$
Because $I_z+F_z=\overline{\bigcup_{k<0}z^kI_z}$,
we'll have checked the Theorem \ref{main} then.

\subsection{The general behavior}
Fix some $n,D,z$ as in Observation \ref{sps1}.\\
We have $n\geq 2,\ D^2<4n$ - it's easy to check: \be \label{ln}
D\leq n,\quad n=D\Rightarrow n\in\{2,3\} \ee

The equation $z^2=Dz-n$ can be written as \be\label{zero}(n,-D,1)_z=0 \ee where we use
the natural identification: $$\mathbb{Z}^k\ni (b_0,..b_{k-1})\equiv
(a_i):\forall_{i<0\vee i\geq k}\ a_i=0,\ \forall_{0\leq i
<k}\ a_i=b_i.$$
Now take some $\mathbb{C}\ni x:=(b_0,...,b_k)_z$, where $\forall_i b_i\in\mathbb{Z}$.\\
If we add to this pseudorepresentation of $x$ any multiplicity of $(n,-D,1)$ or some of its shift,
it will still represent $x$.
If we subtract $\lfloor b_0/n \rfloor (n,-D,1)$, we will make the zeroth digit to
be in $\bar n$.\\
It suggest an algorithm to change a psudorepresentation into
representation - use above construction on succeeding digits to
bring them into $\bar n$.\\
The question is: will it finally stop?
\begin{ex} \emph{We can use $(n,-D,1)_z=0$ in many ways.}
  $$n=3,D=3\Rightarrow z=\frac{3}{2}+j\frac{\sqrt 3}{2},\
  (3,-3,1)_z=0$$
  \emph{We can change some pseudorepresentation into
  representation:}
  $$(7,-8,7,-2)\rightarrow(1,-2,5,-2)\rightarrow(1,1,2,-1)\rightarrow
  ...\rightarrow(1,1,0,...,0,2,-1)$$
  \emph{Find the representation of some $x\in X$:}
  $$x=10z-12=(-12,10)_z=(0,-2,4)_z=(0,1,1,1)_z=z+z^2+z^3$$
  \emph{Or find the position on $X$ for some representation:}
  $$z^4=(0,0,0,0,1)_z=(0,0,-3,3)_z=(0,-9,6)_z=(-18,9)_z=9z-18,$$
  \emph{generally:} $z^k=\left( \begin{array}{cc} z& 1 \end{array}\right)
\left(\begin{array}{cc}D & 1 \\-n & 0 \\\end{array}\right)^k
\left(\begin{array}{c}0 \\1 \\\end{array}\right).
$
\end{ex}

 \textbf{Digression} We can use arithmetic in such bases similar as usual -
 for example adding would be adding two representations position by position -
 we get a pseudorepresentation of the sum. We have to change it into representation.
 In the binary system we would use $(2,-1)_2=0$, this time $(n,-D,1)_z=0$.\\
 We can create parallel algorithm for it - in the binary system we
 divide the
 representation into two parts and sum them parallelly, for every possible carried digit (0 or 1).
 The difference in this case, is that we have more 'carried cases':
 for addition from 8 (for $D=1$) to 21 (for $n=3,\ D=3$).\\

$x\in \mathbb{Z}+z\mathbb{Z}=x_0+zx_1$ can be 'psuedorepresented' as
$x=(x_0,x_1)_z$, so
\begin{lem}
  \emph{For any }$x\in X$, \emph{there exists its} reduction $r(x):=y\in X$, \emph{such that:}
$$x-yz\in \bar{n}.$$
\end{lem}
We can think about the reduction as throwing out the youngest
digit($\in\bar n$).\\

We will now check if $r$ have a nonzero fixed point(s).\\
Assume that $r(x)=x$ for some $x=a+zb$:\\
$$r\left((a,b,0)_z\right)=r\left((a-kn,b+kD,-k)_z\right)=
(b+kD,-k)_z$$
where $k=\lfloor a/n\rfloor$.\\
$a+bz=(b+kD)-kz$, so $b=-k$, $a$ and $b$ has different signs,\\
$a=b+kD=b-bD=b(1-D)$ - we must have $D>1$,\\
$b=-\lfloor b(1-D)/n\rfloor=\lceil b\frac{D-1}{n}\rceil$, but
$0<\frac{D-1}{n}<1$ (\ref{ln}) so $b>0$ \\
for $b=1$ we have $(1-D,1)$ for any $D\geq2$,\\
for $b=2$, we would need $2(D-1)>n$ - it can happen only for $n=3,\
D=3$ - it's very special case we would focus later on.\\
Because of (\ref{ln}) it's the largest possible $b$.\\

We have that there are 1-3 fixed points of $r$. We would like to
show that from all points of $X$ we would get in finite number of
reductions to one of these points.
\begin{lem} \label{red}
  $|x|>\sqrt n+1 \Rightarrow |r(x)|<|x|$
\end{lem}
\textbf{Proof:} $|r(x)|=\left| \frac{x-a}{z}\right|\leq
\frac{|x|+(n-1)}{|z|}=\frac{|x|+(\sqrt n +1)(\sqrt n -1)}
{\sqrt n}<\frac{|x|+|x|(\sqrt n-1)}{\sqrt n}=|x|.\quad\quad \Box$\\

Because we are in the discrete lattice $X$, starting from any point
$x\in X$ after a finite number of reductions, we will get into the ball
$B:=\{x\in X: |x|<\sqrt n +1\}$.\\
It's finite set, all fixed points are in it - if we would show
that there are no cycles in it, the reduction process will always
stabilize.
\begin{df}
  Attractor \emph{of fixed point $s$, is}
  $A_x:=\{x:\exists_i  r^i(x)=s\}$.
\end{df}

\begin{figure}[h]
    \centering
        \includegraphics[width=15cm]{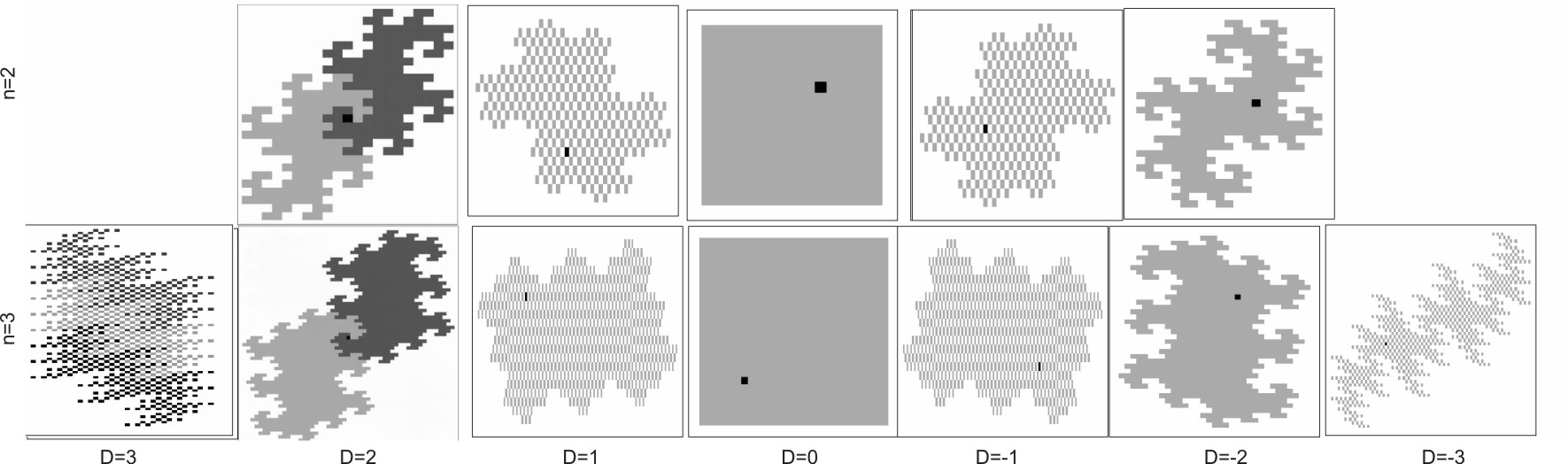}
        \caption{All cases for $n=2,3$. Different colors represents different attractors.
        Black dot is the center of coordinates.} \label{ex2}
\end{figure}

$B$ is finite, so for a finite number of cases we can do it by
checking that in fact all points of $B$ will be in some attractor.
I've done it for $n=2,3$ (Fig. \ref{ex2}).

For $n\geq 4$, $\sqrt n+1\leq n-1$, so $\bar{n}+z\bar{n}$, which is
in the attractor of $0$, covers $\{x\in B: \textrm{arg}(x)\in
[0,\varphi]\}$ (see Fig. 5).\\
$\varphi\in (0,\pi)$ : if $\textrm{arg}(x)\in[\pi,2\pi)$,
$\textrm{arg}(r(x))<\textrm{arg}(x)$.\\
So finally this cycle could be only in $\{x\in B:\textrm{arg}(x)\in
(\varphi,\pi)\}.$\\
If $\varphi>\pi/2,\ (D<0)$, after one reduction we get into the
$\bar{n}+z\bar{n}$, else after 2 reductions we get into the second
stationary point.\\

\begin{figure}[h]
    \centering
        \includegraphics[width=15cm]{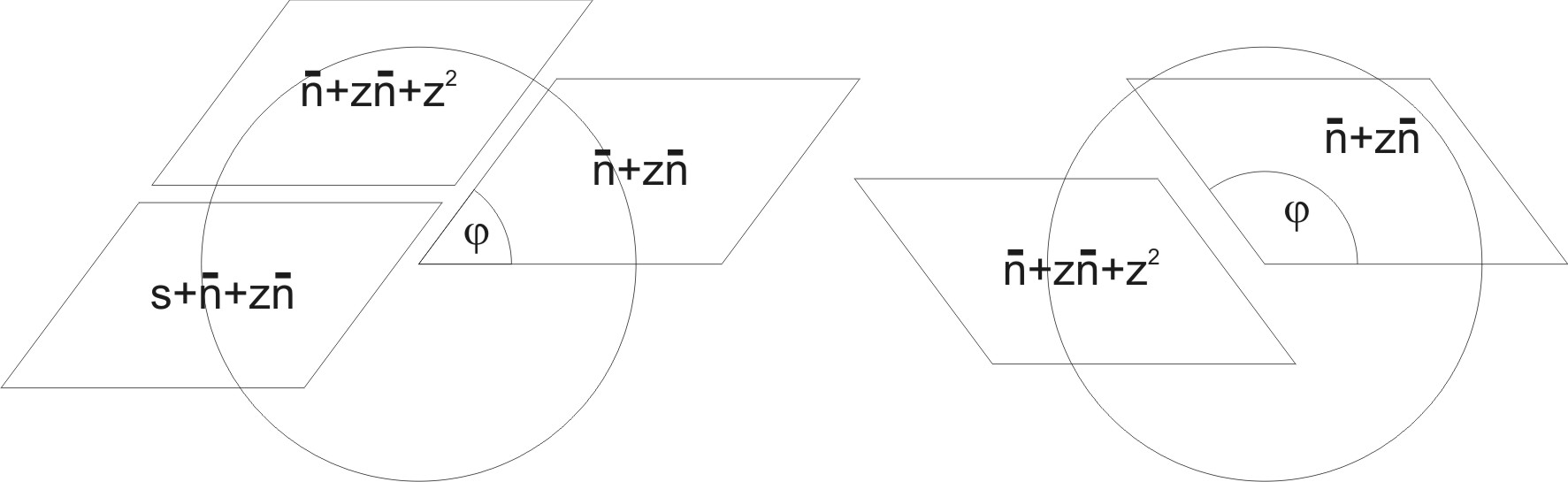}
        \caption{Constructing a cycle.}
\end{figure}
We've just proven Theorem \ref{main}.  $\quad\quad\Box$

\subsection{Special cases}
We will now concentrate on the $D\geq 2$ cases: when $I_z\neq X$.\\

In this case we can split $X$ into two or three disjoined subsets (attractors of fixed points) -
we can still use is as a numeral system, but this time we
have to remember the number of fixed point too
(or of cycle for some cases in higher dimensions).\\
For each fixed point $s=r(s)$, we can split the attractor into subsets ($A_k$)
of points which have the same number of reduction ($k$) to get to the fixed point
(we need $k$ digits to encode it):
$$zA_s^k+\bar{n}=A_s^{k+1},\quad r(A_s^{k+1})=r(A_s)$$
So because $A_s^0=\{s\}$,
$$A_s^k:=sz^k+\sum_{0\leq i<k} z^l\bar{n}\quad(=sz^k+A_0^k)$$
Using $A_s=\sum_{k\geq 0} A^k_s\quad A_s^0\subset
A_s^1\subset...\subset A_s^k\subset...$, we have
\be zA_s+\bar{n}=A_s \ee
it's the same similarity equation as for $I_z$ (\ref{ss2}).\\

There can be some nice symmetries here: check for which $x_0$,
$x_0-I_z$ fulfills (\ref{ss2})
$$z(x_0-I_z)+\bar{n}=z(x_0-I_z)+n-1-\bar{n}=zx_0+n-1-(zI_z+\bar{n})=zx_0+n-1-I_z$$
If we want the right side to be equal $x_0-I_z$, we must have
$zx_0+n-1=x_0$ \be x_0=\frac{n-1}{1-z} \ee For $D=2$, we have
$x_0=\frac{n-1}{1-D-j\sqrt{n-1}}=j\sqrt{n-1}=z-1\in X$,\\
For $n=3,\ D=3$, $x_0=2z-4$ we can say, that $A_{z-2}$ separates this two copies.\\

\begin{figure}[h]
    \centering
        \includegraphics[width=15cm]{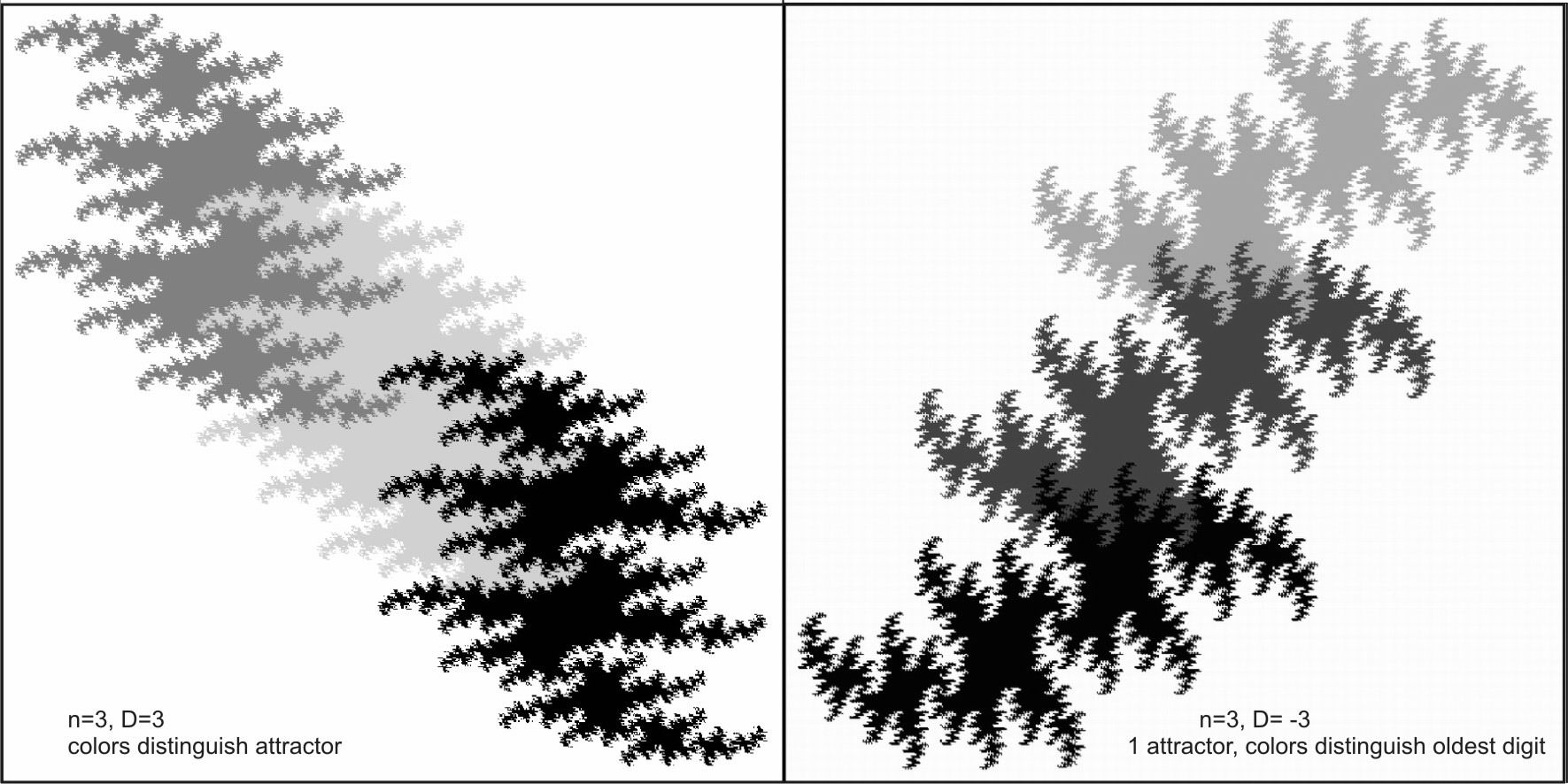}
        \caption{The difference for $D = \pm 3$.}
\end{figure}

Unfortunately for the rest cases $x_0\notin X$ - we can see it using
$$A_0^k=\sum_{i=0..k-1}z^i\bar{n}=\sum_{i=0..k-1}z^i(n-1-\bar{n})=\sum_{i=0..k-1}z^i
(n-1)-A_0^k=(n-1)\frac{z^k-1}{z-1}-A_0^k$$ so
$A_s^k=sz^k+(n-1)\frac{z^k-1}{z-1}-A_0^k$, but this point moves with
$k$ in the other cases.

Finally
\begin{tw} \emph{We can have the following periodic cases:}
\begin{itemize}
  \item \emph{For} $D\leq 1$, $X=A_0$,
  \item \emph{for} $D=2$, $X=A_0\cup A_{z-D+1}=I_z\cup (z-1-I_z)$,
  \item \emph{for} $n=3,\ D=3$, $X=A_0\cup A_{z-2}\cup A_{2z-4}=I_z\cup
  (2z-4-I_z)\cup A_{z-2}$
  \item \emph{for the rest of cases} $X=A_0\cup A_{z-D+1}$.
\end{itemize}
\end{tw}

\section{The fractional part}
In this section will be shown the methodology to analyze $F_z$ - attractors of simple
iterated functional systems.

\subsection{The convex hull}
I will shortly present the results from \cite{me}, which allows to find analytically
the convex hull of such simple fractals.\\
The idea is to define some function which can be easily rewritten in
the selfsimilarity form -
we get some functional equations, which can be solved analytically or approximated.\\
The function we need in this case, is the width function, which gives for every direction
the width in that direction.\\
The other functions which can be threaded in such way are for example:
$f(x)=\mu(F_z\cap (F_z+x))$ (some Husdorff's measure now),\\
$f_p(x)=\sum_{y\in I_z} (x-y)^{-p}$ (for example $p=2$)\\
$\tilde{f}(x)=\sum_{y\in I_z} e^{jx\cdot y}(\textrm{or}\ =\int_{F_z}e^{jx\cdot y}dy)$
(this time '$\cdot$' means scalar multiplication). \\

In our cases, the equations for the width function can be solved analytically.\\
They defines the convex hull of a set as the intersection of all halfplanes:
$$\textrm{conv}(F_z-x_0)=\bigcap_{\alpha\in [0,2\pi)} \{x:\Re\left(e^{-j\alpha}x\right)\leq h(\alpha)\}$$
$$h(\alpha)=\frac{n-1}{2}\sum_{j>0}r^{-j}|\cos(\alpha+j\varphi)|$$
$$x_0=\frac{1}{2}\frac{n-1}{z-1}$$
where $x_0$ is its center of symmetry, $h(\alpha)$ gives the
position of bounding line in $\alpha$ direction. There is also shown
how to construct analytically this convex hull from triangles, find
the length of its boundary $\left( 2\sqrt n +2\right)$ or its area
$\left((n-1)\sum_{i>0}|\sin(i\varphi)|\sqrt n^{-i}\right)$.

\subsection{The construction of the boundary}

In this subsection, will be shown the methodology of using
succeeding approximations of $F_z$:
$$F_z=\frac{\sum_{i=1,..,k}z^k\bar{n}+F_z}{z^k}=
\overline{\bigcup_{k>0}\frac{\sum_{i=1,..,k}z^k\bar{n}}{z^k}}$$ the
first form is true for any $k\geq 0$ - we can think about $F_z$ as a
discrete net of smaller copies of $F_z$. \\
The second: sum of growing
family of such nets, shows that taking the limit $k\to \infty$, we
can 'forget' about \emph{tiles} (copies of $F_z$) - approximations
restores the whole set.\\
We will introduce discrete versions of some topological properties,
like path, boundary, connectiveness, show how to make the step to
the next approximation and finally that in the
limit, they really corresponds to their continuous equivalent.\\

We will finally prove:
\begin{tw} \label{frth}
For periodic cases $\left(z^2=Dz-n,\ D\in 2(-\sqrt{n},\sqrt{n})\right)$, we have:
\begin{itemize}
\item if $D\neq 0$, we can split (intersecting in at most one point) $\delta F_z$ into six connected components - common
parts with its 'neighboring tiles',
\item $F_z$ - connected, simply connected,
\item $\overline{\textrm{int}(F_z)}=F_z.$
\end{itemize}
\end{tw}

The purely imaginary case: $z=j\sqrt n\ (D=0)$ can be easily solved
\cite{knu} - even digits corresponds real part, odd: imaginary,
$F_z$ is just a rectangle. In this case $F_z+I_z$ behave not like
for the others: hexagonal lattice, but it's just a rectangular
lattice - we will omit this case in this subsection.\\
Look on (\ref{neg}) - changing the sign of $D$ will only transpose
$F_z$ - in this subsection we can restrict to \be D<0.\ee

For approximations we can work with $X=\mathbb{Z}+z\mathbb{Z}$ space:
\begin{df} \emph{For any $A\subset X$ define:}\\
  Magnification \emph{of} $A$ \emph{is} $M(A)=Mz+\bar{n}=\bigcup_{i=0,..,n-1}(Mz+i)$.\\
  k-th approximation \emph{of $A$} is $M^k(A)/z^k$.
\end{df}
We can think about it that we see only points from discrete lattice ($X$) and
we use magnifying glass, which allows us to increase magnification
$|z|$ times in one step (plus rotation).\\
After such step, every point(tile) occurs to be $n$ points(tiles).\\

\begin{figure}[h]
    \centering
        \includegraphics[width=15cm]{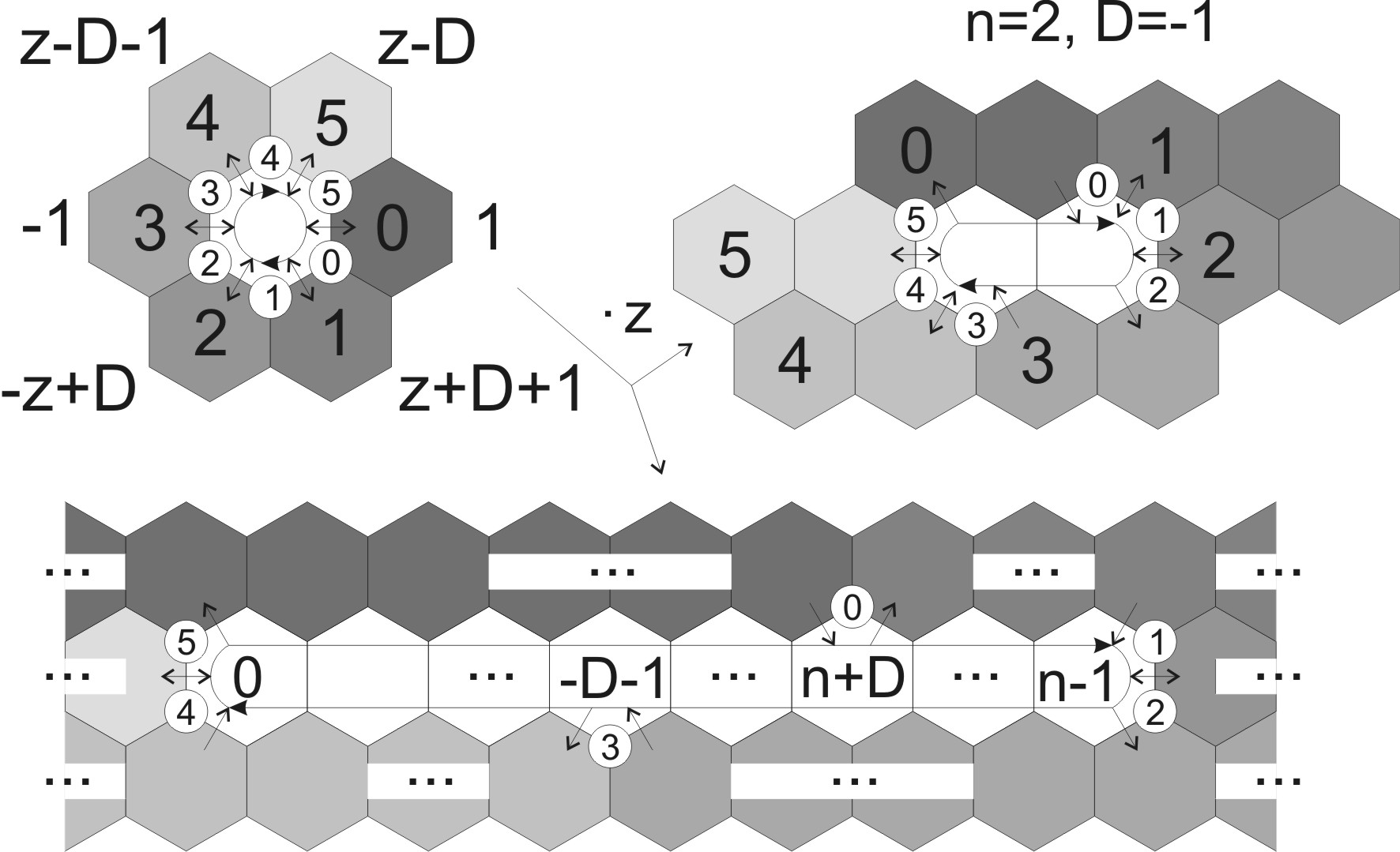}
        \caption{The construction of magnification.} \label{sch}
\end{figure}

To have the hexagonal behavior, we would need to describe neighbors
(Fig. \ref{sch}).\\
I will just give positions and in the next lemma will be shown how
to check them.

\begin{df}
\ \\
Neighbors index set: $R=\mathbb{Z}_6\quad (\mathbb{Z} \textrm{\emph{ modulo }} 6)$.\\
$r$-th neighbor \emph{of $x\in X$ (where $r\in R$) is $N_r(x):=x+N_r$, where:}\\
$N_0:=0,N_1:=-z+D+1,N_2:=-z+D, N_3:=-1, N_4:=z-D-1, N_5:=z-D.$
\emph{The} $r$-th meeting place \emph{is $C_r(x):=xz+C_r$, where}\\
$C_0:=n+d,C_1:=n-1,C_2:=n-1,C_3:=-D-1,C_4:=0,C_5:=0$.\\
\emph{The} $r$-th meeting point \emph{is $P_r(x):=x+P_r$, where}\\
$P_0:=1+\frac{n}{z-1},P_1:=1+\frac{n-D}{z-1},P_2:=\frac{n}{z-1},P_3:=\frac{1}{z-1}-1,
P_4:=\frac{D-1}{z-1}-1,P_5:=\frac{1}{z-1}$.
\end{df}
$N$ are the positions of centers of neighboring tiles - we order
them clockwise. After one step of magnification, the center tile is
split into $n$ tiles, $C_r$ tells in which of them is the common
point with the magnifications of $r$-th and $r+1$-th neighbor (dots with numbers
on the picture).\\
We can easily find the real point of this meeting - do the
magnification:
$$zP_r=P_r+\textrm{ the position of tile for which $zP_r$ is $r$-th place of meeting}$$
for example $zP_0=P_0+(n+D)+(z-D-1)$.

We will now check properties for making a step of magnification.
\begin{lem} For $x\in X,\ r\in R$, we
have that $C_r(x)\in M(x)$ is the neighbor of some elements from
$M(N_r(x))$ and $M(N_{r+1}(x))$.
\end{lem}
\emph{Proof:} We will check it for $x=0$ - we can translate everything to get the thesis.\\
It's easy to check the table below ($z^2=Dz-n$):

\begin{tabular}{|c|c|c|c|c|c|c|}
  \hline
  $r$ & 0 & 1 & 2 & 3 & 4 & 5 \\ \hline
  $N_r$ & 1 & $-z+D+1$ & $-z+D$ & -1 & $z-D-1$ & $z-D$ \\
  $zN_r$ & $z$ & $n+z$ & $n$ & $-z$ & $-n-z$ & $-n$ \\ \hline
  $C_r$ & $n+D$ & $n-1$ & $n-1$ & $-D-1$ & 0 & 0 \\
  $a$ & $n-1$ & $-D-1$ & 0 & 0 & $n+D$ & 0 \\
  $r_a$ & 4 & 5 & 0 & 1 & 2 & 3 \\ \hline
  $b$ & 0 & 0 & $n+D$ & $n-1$ & 0 & $-D-1$ \\
  $r_b$ & 5 & 0 & 1 & 2 & 3 & 4 \\ \hline
  $P_r$ & $1+\frac{n}{z-1}$ & $1+\frac{n-D}{z-1}$ & $\frac{n}{z-1}$
 & $\frac{1}{z-1}-1$ & $\frac{D-1}{z-1}-1$ & $-\frac{1}{z-1}$ \\
  \hline
\end{tabular}

\noindent such that $N_{r_a}=zN_r+a,\ N_{r_b}=zN_{r+1}+b$,\\
for example for $r=0:\ N_{r_a}(C_r)=n+D+N_4=z+(n-1)=zN_r+a\quad \Box$.\\

\begin{df}
\ \\
\emph{The} neighborhood \emph{of $x\in X$ is $N(x):=\{N_r(x):r\in R\}$}.\\
\emph{The} edge \emph{of $A\subset X$ is $\left(\bigcup_{x\in A}
N(x)\right)\backslash A$}.\\
\emph{A sequence $(x_i)_{0\leq i<k}\ (x_i\in B\subset x)$ will be
called a }length $k$ chain in $B$\emph{ if succeeding elements
are neighbors.}\\
Chain \emph{is called } closed\emph{, when additionally $x_0$ and
$x_{k-1}$ are neighbors.}\\
Subchain \emph{of some closed, length $k$ chain if is created from
some of its succeeding elements (0 is after $k-1$).}\\
\emph{A set $A\subset X$ is called} connected\emph{, if any two of its
points can be connected by a chain in this set.}\\
\emph{A set $A\subset X$ is called} simply connected\emph{, if
$X\backslash A$ is connected.}\\
The range of indexes: $[a,b]$\emph{, where $a,b\in R$, is called a
step 1 arithmetic sequence in $R$ from $a$ to $b$ (eg.
$[4,2]=(4,5,0,1,2)$ ).}\\
\emph{A set $B\subset X$} lies on the left of \emph{sequence
$(x_i)_{0\leq i<k}$, when for each $1\leq i\leq k-2$:\\
if $x_{i-1}$ and $x_{i+1}$ are correspondingly $a$ and $b$ neighbor
of $x_i$, then}
$$\{x_i+N_k:k\in[a+1,b-1]\}\subset B.$$
\emph{A set $B\subset X$} lies on the left of \emph{closed chain, if
this property is fulfilled for each of its subchains.}\\
\end{df}

We need lying on the left to make closed chain around some set oriented (counterclockwise).\\
The approximation of the boundary of $F_z$ can be started from
length 6 closed chain - neighbors of 0 in counterclockwise order -
$\{0\}$ lies on the left of this chain.\\
Now we have to make a general construction of making the step to the
next magnification.
\begin{lem} \label{app}
 \emph{Let $B\subset X$ lies on the left of some (closed) chain
 $(x_i)_{0\leq i<k}$.\\
 Then there exists } (closed) chain $(y_i)_{0\leq i<l}$, \emph{such that
 $y_0\in M(x_0),\ y_{l-1}\in M(x_{k-1})$ and that $P(B)$ lies on the
 left of it.}
\end{lem}
\emph{Proof:} Firstly we will find $(y_i)$ for length tree chain:
$(N_a,0,N_b)$ and set $\{N_{a+1},..,N_{b-1}\}$.\\
Define $y_0\in M(N_a)$ to be neighbor of $C_a$, $y_{l-1}\in
M(N_b)$ to be neighbor of $C_b$.\\
The rest of the chain lies in $M(0)=\bar{n}$, such that the set lies
on the left - it's chosen unambiguously - look on Fig. \ref{sch} :
using some arrow we enter $\bar{n}$ in $C_a$, go along clockwise cycle
inside, to exit turning left in $C_b$.\\
For example (denote $(N_a,0,N_b)$ as $a_b$): $0_4\to
4_0(3_0)^{-D-2}3_3 (0_3)^{n+D-1}0_2$.\\
In the general situation, we are taking succeeding $x_i$ and use the
above construction for $(x_{i-1}-x_i,0,x_{i+1}-x_i$ (remembering to
remove duplicates on the ends).  $\Box$\\

Now starting from $C_0$ - closed, counterclockwise chain around
$\{0\}$, call $C_i$ - the $i$-th use of above lemma to $C_0$.\\
It's the edge of $M^i(\{0\})$, which lies on the left of
$C_i$.\\
When we take approximations (divide this situation by $z^i$), we see
that we are getting closer to $F_z$ and its (oriented) boundary.\\
$$C=\{x\in \mathbb{C}: \exists_{x_i}\ x_i\to x,\ x_i\in C_i/z^i\}$$
 \emph{Proof of Theorem \ref{frth} :} We have only to check:
\begin{itemize}
  \item $C$ is closed - just take diagonal sequence.
  \item $C$ is connected - the distance between
  succeeding points of chain $C_i/z^i$ goes to zero.
  \item $\mathbb{C}\backslash C$ has exactly two connected
  components - it's true for every approximation, using the previous
  point we have it in the limit.
  \item The interior component with $C$ is $F_z$, $C=\delta F_z$ -
  $M^i(\{0\})/z^i$ are approximations of $F_z$ and its distance to
  $C_i/z^i$ is going to 0.
  \item In each step we can divide $C_i/z^i$ into six subchains -
  lying in approximations of succeeding neighbors ($M^i(N_a)/z^i$) -
  we can split $C$ into 6 subsets this way.
  \item $\overline{\textrm{int}(F_z)}=F_z$ - from the Bair's theorem: if int$(z)=\emptyset$,
then the sum of countable number of them couldn't give the whole space.\\
Because $\forall_k\ F_z=(M^k(\{0\})+F_z)/z^k$, we can cover $F_z$ with its smaller copies
having arbitrary small diameter - some of them will lie in this nonempty interior.
After rescaling back, we've got the thesis.
$\quad\quad\quad\quad\quad\quad\Box$
\end{itemize}

Using the above construction of succeeding approximations of $\delta
F_z$, we can for example calculate its Hausdorff's dimension.\\
Namely using Lemma \ref{app}, we can change two edges: $(N_a,0,N_b)$
into a sequence of them in the next magnification.\\
If for each edge we distinguish between the direction of the next
edge - remember the first of this two edges as $a_b$, we can assign to
it a sequence of such pairs in the next magnification
(like in the example, but without the last $0_2$).\\
If we are interested in the number of edges only (instead of above grammar),
we can write this
iteration in the $36\times 36$ (or $18\times 18$ if we use symmetry)
matrix form. It's normalized dominant eigenvector tells us the
asymptotic probability distribution of edges, corresponding
eigenvalue($\lambda$) tells that asymptotically, for the next
approximation we need $\lambda$ times more edges, but the length of
them is $\sqrt n$ times smaller, so the boundary Hausdorff's
dimension is:
$$H=\frac{\log(\lambda)}{\log(\sqrt n)}$$
Here are found values for $n\leq 9$:

\begin{tabular}{|c|c|c|c|c|c|c|c|c|c|}
  \hline
  $D/$ & 0 & 1 & 2 & 3 & 4 & 5 \\
  $n$ &    &   &   &   &   &      \\ \hline
  2   & 1 & 1.210760533 & 1.523627086 &  &  &      \\ \hline
  3   & 1 & 1.162039854 & 1.376841713 & 1.657559542 &  &      \\ \hline
  4   & 1 & 1.134761994 & 1.303052340 & 1.508664987 &  &      \\ \hline
  5   & 1 & 1.116924317 & 1.257583258 & 1.422944863 & 1.608726378 &      \\ \hline
  6   & 1 & 1.104171451 & 1.226294386 & 1.366294523 & 1.520716574 &      \\ \hline
  7   & 1 & 1.094508825 & 1.203216411 & 1.325629733 & 1.458928649 & 1.598134771     \\ \hline
  8   & 1 & 1.086882303 & 1.185363533 & 1.294784127 & 1.412801706 & 1.535582008     \\ \hline
  9   & 1 & 1.080677473 & 1.171064018 & 1.270444950 & 1.376841713 & 1.487192945     \\
  \hline
\end{tabular}
\ \\
In our cases, $\delta F_z$ in fact - defines $F_z$, and it's
dimension is smaller.\\
So we can use just the boundary - operating on it (drawing for
example) is much faster.

\section{Higher dimensions}
In this informal section there will be shortly introduce a suggestion
of higher dimensions generalization, focusing on the dimension 3.\\
Once more I will do it only for the periodic case. \\
I haven't even checked the hypothesis numerically in this case.\\

We can think about $z^k$ from the previous sections as
\be \label{dff}z^k:=Z^k\mathbf{1}\quad\quad\quad Z^k\mathbf{1}\cdot Z^l\mathbf{1}:=
(z^kz^l=z^{k+l}=)Z^{k+l}\mathbf{1} \ee
where $\mathbf{1}=(1,0)^T,\quad Z=r\left(\begin{array}{cc}\cos(\varphi) & -\sin(\varphi) \\
\sin(\varphi) & \cos(\varphi) \\ \end{array} \right)$.\\
We can use this definition (\ref{dff}) in higher dimension Euclidean space.\\
We will focus on $\mathbb{R}^3$.\\

Choose some unit vector ($\mathbf{1}$), $r\in \mathbb{R}:\ |r|>1$ and
some orthogonal matrix $O\in \mathbb{R}^{3\times 3}$.\\
Define $$Z:=rO$$
similarly as for complex numbers.\\
We can diagonalize $O$ - it is a rotation ($\varphi$ radians) around some vector ($u$) and
eventually the symmetry, but we can put it in the sign of $r$.\\
We want $1,z,z^2$ to generate whole space, so the angle between $\mathbf{1}$ and $u$
(denote it $\psi$) cannot be the multiplicity of $\pi/2$.\\
So we can define orthonormal basis of our space: $e_1:=\mathbf{1}$, $e_2$ is orthogonal to $e_1$
and $u$ lies on the space generated by them, $e_3$ is orthogonal to $e_1$ and $e_2$.\\
Define the another basis $\{u',v,w\}$:\\
$u'=\cos(\psi)\cdot\left(\cos(\psi),\sin(\psi),0\right)$ - rescalled $u$,\\
$v=\sin(\psi)\cdot\left(\sin(\psi),-\cos(\psi),0\right)$ - orthogonal to $u$,\\
$w=e_3$.\\
Now we can join $v$ and $w$ and think about our space as $\mathbb{R}\times\mathbb{C}$:\\
$\mathbf{1}=(1,1)$, we have multiplication and addition by
coordinates, \be z^k=(r^k,r^ke^{jk\varphi}) \ee We can define the
analog of the conjugancy: $\{e_1,e_2\}$ plane symmetry
$\left(\overline{(x,y)}:=(x,\bar{y})\right)$.\\

Analogically to (\ref{powl}), we have the condition:
$$r=\sqrt[N]{n}$$
in $N$ dimensional space.\\

To assure the periodicity in dimension 3, we would need ($|r|^3=n$):
\be z^3=-Az^2-Bz-C\mathbf{1} \ee
for some $A,B,C\in \mathbb{Z}$.\\
conjugate this equation and multiply by $z^3$:
$$r^6\mathbf{1}=-Ar^4z-Br^2z^2-Cz^3$$
comparing both equations, we get (the sign of $C$ is hidden in $r$):
$$C=r^6/C, \quad A=Br^2/C,\quad B=Ar^4/C$$
$$C=-r^3,\quad B=-Ar$$
We see that $r$ have to be integer. We've found the analog of (\ref{zero}):
\be (-r^3,-Ar,A,1)_z=0 \label{zeroo}\ee

We now have to find some $z=(R,Re^{j\varphi})$ fulfilling this equation.\\
The real part gives:
$$0=-r^3-ArR+AR^2+R^3$$
$R\in \mathbb{R}$, so $R=r$.\\
The complex part gives now:
$$0=-r^3-Ar^2e^{j\varphi}+Ar^2e^{2j\varphi}+r^3e^{3j\varphi}=
r^2\left(r(e^{3j\varphi}-1)+Ae^{j\varphi}(e^{j\varphi}-1)\right)$$
$$-A/r=e^{-j\varphi}(e^{2j\varphi}+e^{j\varphi}+1)=
e^{j\varphi}+1+e^{-j\varphi}=1+2\cos(\varphi)$$
\be \varphi=\arccos\left(-\frac{1}{2}(1+A/r)\right)\quad\quad
\left(\Rightarrow A\in\{-3r+1,...,r-1\}\right)\ee
The sign of $\varphi$ is only the matter of conjugation.\\

Using (\ref{zeroo}) we can make sequence of reductions like in 3.1.\\
By analogy to Lemma \ref{red}, we have
\be|x|>\frac{n-1}{\sqrt[N]n-1} \Rightarrow |r(x)|<|x|\ee we will
finally get into this ball, and finally
get to a fixed point or a cycle.\\
We can check that $(0,0,0)$ is the only fixed point for $A\geq 0$.\\
But in this case, for $r>0$, we get a length two cycle:\\
$(2A-1,-A.-1)\to (-A,-1,0)\to...$.\\
I've checked, that in the second case:
\be m\geq 2,\ r=-m,\ n=m^3,\
A\in\{0,...,3m-1\},\quad (n,mA,A,1)_z=0\ee for a few first $m$ everything
is fine - we can generate the whole ball starting from $(0,0,0)$ -
we have proper numeral systems.\\

We can generalize the methods for the fractional part in this cases
too.\\
The width function in introduced coordinates have two arguments, say
$\alpha \in[-\pi/2,\pi/2],\ \beta\in[0,2\pi)$ which corresponds to
$(\cos(\alpha),\sin(\alpha)e^{j\beta})$.\\
Now using $z^2F_z=F_z+\bar{n}(1,1)+\bar{n}(r,re^{j\varphi})$, we
have equations for $h$, in which $\alpha$ is fixed - we can find
solutions as a infinite sum for a fixed $\alpha$, like in \cite{me}.\\
They don't correspond to one plane (cross section) this time - to
find $u$ coordinate, we need to use $dh/d\alpha$.\\

We can use the method to construct the boundary of $F_z$.\\
This time every tile has 14 neighbors:\\ $N(0)=\pm \{1,A+z,A-1+z,
Am+Az+z^2,Am-1+Az+z^2,$\\
 \hbox{ } $\quad\quad\quad\quad,Am-A+(A-1)z+z^2,Am-A+1+(A-1)z+z^2\}$ \\
 Now the simplest structure we can separately magnify, analogically
 to
$(N_a,0,N_b)$, is 0 with closed chain made of its neighbors, which
splits $N(0)$ into two connected subsets, which can be distinguished using
chain orientation (one of them can be empty).

\end{document}